\documentclass[10pt,amssymbols]{amsart}

\begin{document}
\bibliographystyle{amsplain}

\title{A correction to: Rational and nearly rational varieties }
\author{J\'anos Koll\'ar, Karen E.\ Smith and Alessio Corti }

\maketitle

In accordance with  a request of A.\ V.\ Pukhlikov,
we  acknowledge that Lemma 5.36 in our book
{\it Rational and nearly rational varieties}
should have been attributed to 
him. The statement and the proof of   Lemma 5.36 was first published in:

A.\ V.\  Pukhlikov: 
{\it A remark on the theorem of V.\ A.\ Iskovskikh and Yu.\ I.\ Manin
            on a three-dimensional quartic},
Trudy Mat. Inst. Steklov. 208 (1995), 278--289

We regret that we treated it as part of mathematical ``folklore''
when a proper reference was available.

\end{document}